\documentclass[12pt,a4paper]{amsart}
\usepackage{amsmath,amssymb,amsfonts,amsthm}

\newtheorem{lemma}{Lemma}
\newtheorem{remark}{Remark}
\newtheorem{condition}{Condition}
\newtheorem{theorem}{Theorem}
\newtheorem{prop}{Proposition}

\newcommand{\cs}{^*}

\newcommand{\biz}{\textbf{Proof.\ }}
\newcommand{\ri}{\right)}
\newcommand{\lef}{\left(}

\newcommand{\lk}{\left\lbrace}
\newcommand{\rk}{\right\rbrace}

\newcommand{\Prob}{\mathbb{P}}
\newcommand{\e}{\varepsilon}
\newcommand{\lsz}{\left[}
\newcommand{\rsz}{\right]}

\newcommand{\F}{\mathcal{F}}
\newcommand{\tv}{\rightarrow \infty}
\newcommand{\aaa}{\textit a$)$}
\newcommand{\bb}{\textit b$)$}
\newcommand{\al}{\alpha}

\title{Local Degree Distribution in Scale Free Random Graphs}
 
\author{\'Agnes Backhausz}
\address{Department of Probability Theory and Statistics\\Faculty of Science\\
 E\"otv\"os Lor\'and University\\P\'azm\'any P.~s.\ 1/C, H-1117
 Budapest, Hungary} 
\email{agnes@cs.elte.hu}
\author{Tam\'as F.~M\'ori}
\address{Department of Probability Theory and Statistics\\Faculty of Science\\
E\"otv\"os Lor\'and University\\P\'azm\'any P.~s. 1/C, 
H-1117 Budapest, Hungary}
\email{moritamas@ludens.elte.hu}
\dedicatory{\upshape Department of Probability Theory and Statistics,
Faculty of Science\\ E\"otv\"os Lor\'and University\\ P\'azm\'any P.~s. 1/C,
H-1117 Budapest, Hungary\\                         
\textit{E-mail address:} \texttt{agnes@cs.elte.hu, moritamas@ludens.elte.hu} }

\keywords{Scale free, random graphs, recursive trees, martingales, regular variation} 
 
\subjclass[2000]{60G42, 05C80}

\thanks{Research supported by the Hungarian National Foundation for
Scientific Research, Grant No. K67961}
  
\date{14 July 2010}
\begin{document}

\begin{abstract}
In several scale free graph models the asymptotic degree
  distribution and the characteristic exponent change when only a
  smaller set of vertices is considered. Looking at the common
  properties of these models, we present sufficient conditions for the
  almost sure existence of asymptotic degree distribution constrained
  on the set of selected vertices, and identify the chararteristic
  exponent belonging to it.
\end{abstract}
\begin{titlepage}
\maketitle
\thispagestyle{empty}
\end{titlepage}

\section{Introduction}
Since the end of the nineties several complex real world networks and
their random graph models have been investigated \cite{bollobas, Drmota,
  durrett}. Many of them possess the scale free property: the tail of 
the degree distribution decreases polynomially fast, that is, if $c_d$
denotes the proportion of vertices of degree $d$, then $c_d\approx C\cdot d^{-\gamma}$ holds for large
values of $d$ 
\cite{[Bar99]}. $\gamma$ is called the characteristic exponent. 

If the whole network is completely known, the empirical estimator of
the characteristic exponent may have nice properties. However, real
world networks usually are too large and complex, hence our knowledge
of the graph is partial.  For several models of evolving random graphs
the degree distribution and the characteristic exponent change when
attention is restricted to a set of selected vertices that are close
to the initial configuration \cite{[Mo06], [Mo07], [Mo09]}.

Starting from these phenomena, in this paper the degree distribution
constrained on a set of selected vertices will be investigated, assuming
that the graph model possesses the scale free property with
characteristic exponent 
$\gamma>1$, and the number of selected vertices grows regularly
with exponent $0<\alpha\leq 1$. Sufficient conditions for the almost
sure existence of the local asymptotic degree distribution will be
given. It will be shown that under these conditions the
characteristic exponent of the constrained degree distribution is
$\alpha\left(\gamma-1\right)+1$. 

The proofs are based on the methods of martingale theory. Applications
of the general results 
to different graph models (e.g. to the Albert--Barab{\'a}si random
tree \cite{[Bar99]}) will be shown. 

In Section 2 we present the family of random graph models to be
examined and formulate the sufficient conditions. In Sections 3 and 4
we mention some results about martingales and slowly varying sequences
to be applied in the proofs. Section 5 contains the proof of the
main results, and in Section 6 we give some examples and applications.  

\section{Main results}
In this section we present sufficient conditions for the almost sure
existence of asymptotic degree distribution constrained on the set of
selected vertices, and we describe that distribution. 

Let $\lef G_n=\left( V_n, E_n\right)\ri_{n\in\mathbb N}$ be a sequence of
evolving simple random graphs. Some vertices are distinghuised; let
$S_n\subseteq V_n$ denote the set of selected vertices. 

We start from a finite, simple graph $G_0=\lef  V_0, E_0\ri$, this is the
initial configuration with $V_0=\lk u_1, u_2,\ldots,
u_l\rk$. $S_0\subseteq V_0$ is arbitrarily chosen.   
For $n\geq 1$, at the $n\,$th step 
\begin{itemize}
 \item one new vertex, $v_n$, is added to the graph: $V_n=V_0\cup
\lk v_1,\ldots, v_n\rk$; 
\item the new vertex gets some random edges, thus $E_{n-1}\subseteq
  E_n$, and every edge from $E_n\setminus E_{n-1}$ is connected to
  $v_n$; 
\item the new vertex can be added to the set of selected vertices,
  $v_n\in S_n$ is a random choice.  
\end{itemize}

The $\sigma$-field of events generated by the first $n$ steps is
denoted by $\mathcal F_n$. 

For $v\in V_n$ let the degree of $v$ in $G_n$ be denoted by
$\deg_n\left( v\right)$. Furthermore, for $n\geq 1$, and
$d\geq 0$ define 
\begin{align*}
X\left[ n,d\right]&= \left\vert \left\lbrace v\in V_n:
       \deg_n\left(v\right) =d\right\rbrace \right\vert;\\ 
Y \left[ n,d\right]&= \left\vert \left\lbrace v\in V_n: \deg_n\left(
       v\right) =d,\ \left( v, v_{n+1}\right) \in E_{n+1}\right\rbrace
       \right\vert;\\ 
I\left[ n,d\right]&=\begin{cases}
                    1 & \text{if }\deg_n\left( v_n\right)=d,\\
                    0 & \text{otherwise.}          
                        \end{cases}
\end{align*}

In some models it is possible that the new vertex does not get any
edges at some steps. In other models the degree of the new vertex is
fixed, for example, the degree of the new vertex is always $1$ in
random tree models. If the new vertex gets at least $m$ edges at each
step for some $m>0$, then $X\lsz n, d\rsz$ is at most $\left \vert
  V_0\right\vert$ for all $n$ and $d<m$. Thus we denote the minimal
initial degree of the new vertex by $m$, and we consider $X\lsz n,
d\rsz$ only for $d\geq m$. Of course, $m=0$ is also possible. 

\subsection{Conditions on the graph model}

We say that a discrete probability distribution $\left( a_n\right) $
is \emph{exponentially decreasing}  if $a_n\leq C\cdot q^n$ holds for
all $n\geq 1$ for some $C>0$ and $0<q<1$ . A sequence $\left(
  a_n\right)$ is \emph{slowly varying} if
$a_{\left[sn\right]}/a_n\rightarrow 1$ as $n\tv$ for all $s>0$. 

Throughout this paper, for two sequences $\lef a_n\ri, \lef b_n\ri$ of
nonnegative numbers, $a_n\sim b_n$ means that $b_n>0$ except finitely
many terms, and $a_n/b_n\rightarrow 1$ as $n\tv$.  

Now we can formulate the conditions on the graph model.
\begin{condition} \label{c1}
$X\left[ n, d\right]\sim c_d\cdot n$  holds as $n\rightarrow\infty$
for every $d\geq m$ with probability $1$, where $\left( c_d\right)$ is
a probability distribution and $c_d$ is positive for all $d\geq m$. 
\end{condition} 
This means that asymptotic degree distribution exists in this graph
model. Note that $X\left[ n, d\right]\tv$ as $n\tv$ almost surely.

\begin{condition} \label{c2}
$c_d\sim K\cdot d^{-\gamma}$ holds as $d\rightarrow\infty$ for some
positive numbers $K$ and $\gamma$.
\end{condition}
This is the so called scale free property with characteristic exponent
$\gamma$. That is, the asymptotic degree distribution decays
polynomially with exponent $\gamma$. 
This implies that $c_d$ is positive for every $d$ large enough, but we
will need it for all $d\geq m$, this is included in Condition $
\ref{c1}$. 

\begin{condition} \label{c3} 
For every $n\geq 0$, if $w_1, w_2\in V_n$ and $\deg_n\left(
w_1\right)=\deg_n\left( w_2\right)$, then
\[\Prob\left(\left. \left( w_1, v_{n+1}\right)\in E_{n+1}\right\vert
\mathcal F_n\right)= \Prob\left(\left. \left(w_2, v_{n+1}\right)\in
  E_{n+1}\right\vert\mathcal F_n\right).\] 
\end{condition}
In other words, at each step, conditionally on the past, old vertices
of the same degree get connected to the new vertex with the same
probability. 

\begin{condition} \label{c4}
$\sum_{i=1}^{n+1} I\left[ i,d\right] = p_d\cdot n+o\lef n\ri$ holds as
$n\rightarrow\infty$ for every $d\geq m$ with probability $1$, where
$\left(p_d\right)$ is an exponentially decreasing probability
distribution. 
\end{condition}
Loosely speaking, the degree of the new vertex has an exponentially
decreasing asymptotic distribution. This trivially holds if the degree
of the new vertex is fixed. 

\begin{condition} \label{c5}
For every $d\geq m$ there exists a random variable $Z_d\geq 0$ with
exponentially decreasing distribution such that 
\[\Prob \left( Y\left[ i, d\right]\geq l\mid \mathcal F_i\right)\leq
\Prob\left( Z_d\geq l\right), \quad i\geq 1,\ l\geq 1.
\] 
\end{condition}
In many particular cases the following stronger condition is also
met. 

\textit{There exists a random variable $Z\geq 0$ with
exponentially decreasing distribution such that}
\begin{equation}\label{c5strong}
\Prob\lef\deg_n(v_n)\geq l\mid \mathcal F_{n-1}\ri\leq
\Prob\lef Z\geq l\ri, \quad n\geq 1,\ l\geq m.
\end{equation}
This is a sort of upper bound for the initial degree of the new vertex.

\begin{condition} \label{c6}
For every $d\geq m$ we have
\[k_d=\sum_{j=m}^d \lef p_j-c_j\ri>0.\]
\end{condition}
We will see later that the nonnegativity of $k_d$ follows from the
previous conditions; however, the positivity of $k_d$ cannot be
omitted, as an example will show. 

\subsection{Conditions on the set of selected vertices}

Recall that $S_n\subseteq V_n$ is the set of selected vertices in
 $G_n$. We emphasize that  $\deg_n\left( v\right)$ always denotes the
 degree of vertex $v$ in $G_n$, not in $S_n$. 

We will need the following notations. The $\sigma$-field generated by
the first $n$ steps and adding the edges of $v_{n+1}$ at the $\lef
n+1\ri$st step is denoted by $\F_n^+$. Furthermore, for $n\geq 1$ and
$d\geq m$ let       
\begin{align*}
X\cs\left[n,d\right]&= \left\vert \left\lbrace v\in S_n:\ \deg_n\left(
      v\right) =d\right\rbrace \right\vert;\\ 
Y\cs\left[n,d\right]&= \left\vert \left\lbrace v\in S_n:\ \deg_n\left(
      v\right) =d, \left( v, v_{n+1}\right) \in E_{n+1}\right\rbrace
      \right\vert;\\ 
I\cs\left[n,d\right]&=\begin{cases}
                    1 & \text{if $v_n\in S_n$ and $\deg_n\lef v_n\ri=d$,}\\ 
                    0 & \text{otherwise;}          
                        \end{cases}\\
I\cs\left(n \right)&=\sum_{d=m}^n I\cs\left[ n,d\right]=\begin{cases}
                    1 & \text{if $v_n\in S_n$,}\\
                    0 & \text{otherwise.}          
                      \end{cases}
\end{align*}
The conditions on the set of selected vertices are the following.
\begin{condition} \label{d1}
$S_n\subseteq S_{n+1}$ for all $n\geq0$.
\end{condition}
Vertices cannot be deleted from the set of selected vertices.
\begin{condition}\label{d2}
$I\cs\left( n+1\right)$ is $\mathcal F^+_n$-measurable for all $n\geq 0$.
\end{condition}
At each step we have to decide whether the new vertex is to be
selected immediately after choosing its neighbours. Selecting the
neighbours of a fixed vertex is an example. 
 
\begin{condition} \label{d3}
  There exists a sequence of positive random variables
  $\left(\zeta_n\right)$ that are slowly varying as
  $n\rightarrow\infty$, and $\left\vert S_n\right\vert =\sum_{i=1}^n
  I\cs\left(i\right)\sim \zeta_n\cdot n^{\alpha}$ for 
  some $\alpha>0$, with probability $1$. 
\end{condition}
This means that the size of the set of selected vertices is
regularly growing with exponent $\al>0$. 

\begin{condition} \label{d4} For every $d\geq m$ 
\[\sum_{i=1}^n E\left( \left.I\cs\left[i, d\right]\right\vert
\mathcal F_{i-1}\right)= \lef q_d+o\lef 1\ri\ri \sum_{i=1}^n
E\left(\left. I\cs\left(i\right)\right\vert\mathcal F_{i-1}
\right)\] 
holds a.s. as $n\rightarrow\infty$, with some
exponentially decreasing probability distribution
$\left(q_d\right)_{d\ge m}\,$.  
\end{condition}
The last condition holds if the degree of the new vertex $v_{n}$ is
fixed, or its degree and $I\cs\lef n\ri$ are independent conditionally on 
$\mathcal F_{n-1}$. In that case sequence $q_d=p_d$
satisfies the condition. It is also possible that the asymptotic
degree distribution of the new selected vertices is different from
$\left(p_d\right)$ if only it decays exponentially fast. 

\subsection{Description of the local degree distribution}

Now we formulate the main results. 

\begin{theorem}\label{t1}
 Suppose that Conditions \ref{c1}--\ref{d4} hold for a
 random graph model $\lef G_n, S_n\ri$, then the limits 
\[\lim_{n\rightarrow\infty}\frac{X\cs\left[ n, d\right]}{\left\vert S_n\right\vert}=x_d \]
 exist for all $d\geq m$ with probability $1$. 

The constants $x_d$ satisfy the following recursive equations.
\begin{equation*}
x_m=\frac{\al\,q_m}{\al+\tfrac{p_m-c_m}{c_m}}\,,\quad  
x_d=\frac{x_{d-1}\cdot \frac{k_{d-1}}{c_{d-1}}+\alpha\cdot q_d}
{\alpha+\frac{k_d}{c_d}}\quad (d\geq m+1).
\end{equation*}

Sequence $\lef x_d\ri$ is a probability distribution, that is, it sums up
to $1$. Moreover, $x_d\sim L\cdot d^{-\gamma\cs}$  as $d\tv$ with $L>0$ and
\[\gamma\cs=\alpha\left( \gamma-1\right)+1.\]
\end{theorem}

\begin{remark}
From the proof it is clear that with Condition \ref{c2} dropped the
limits $x_d$ still exist and the recursive equations remain valid. 
The role of the scale free property of the graph is just to guarantee that
the asymptotic degree distribution constrained on the set of selected
vertices is also polynomially decaying.  
\end{remark}

\section{Martingales}
We will extensively use the following propositions that are based on
well-known facts of martingale theory. 

\begin{prop}\label{prop0}
 Let $\lef M_n, \mathcal G_n\ri$ be a square integrable martingale with
 $M_1=0$, $\mathcal G_0=\lk \emptyset, \Omega\rk$. Introduce
\[A_n=\sum_{i=2}^n E\lef \left.\lef M_i-M_{i-1}\ri^2\right\vert \mathcal
 G_{i-1}\ri,\] 
that is, the predictable increasing process in the Doob
decomposition of $M_n^2$. Then $M_n=o\lef A_n^{1/2}\log A_n\ri$ holds
almost surely on the event $\lk A_{\infty}=\infty\rk$, and  $M_n$
converges to a finite limit, as $n\rightarrow\infty$, almost surely on
the event $\lk A_{\infty}<\infty\rk$\,. 
\end{prop}

This is a corollary of Propositions VII-2-3 and VII-2-4 of \cite{[Ne75]}.

\begin{prop} \label{prop1} 
Let $\left (M_n, \mathcal G_n\right)$ be a  square integrable
nonnegative submartingale, and   
\[A_n=EM_1+\sum_{i=2}^n \left(E\left(\left. M_i\right\vert\mathcal
    G_{i-1}\right)-M_{i-1}\right), \quad  B_n=\sum_{i=2}^n
\text{\rm{Var}}\left(\left. M_i\right\vert\mathcal G_{i-1}\right).\] 
If $B_n^{1/2}\log B_n=O\left( A_n\right)$, then $M_n\sim A_n$ on the
event $ \left\lbrace A_n\rightarrow\infty\right\rbrace$. 
\end{prop} 

This is easy to prove applying Proposition \ref{prop0} to the
martingale part of the Doob decomposition of $M_n$. 
\begin{prop}\label{prop2} 
Let $Y_1, Y_2,\ldots$ be nonnegative, uniformly bounded random
variables, and $\mathcal G_n=\sigma\lef Y_1, \ldots, Y_n\ri$. Then the
symmetric difference of the events
$\lk\,\sum_{n=1}^{\infty} Y_n<\infty\,\rk$ and $\lk\,
\sum_{n=1}^{\infty} E\lef \left.Y_n\right\vert\mathcal G_{n-1}\ri<\infty
\,\rk$ has probability $0$. Moreover, 
\[\frac{\sum_{n=1}^{\infty} Y_n}{\sum_{n=1}^{\infty} E\lef
  \left.Y_n\right\vert\mathcal G_{n-1}\ri}\rightarrow 1 \quad \lef
n\tv\ri\]
holds almost everywhere on the event $\lk\,\sum_{n=1}^{\infty}
Y_n=\infty\,\rk$.
\end{prop}

This proposition follows from the L\'evy generalization of the Borel--Cantelli
lemma that can be found in \cite{[Ne75]}  (Corollary VII-2-6). 

\section{Slowly varying sequences}

In the proofs we will use the basic results 
of the theory of regularly varying sequences, see e.g.
\cite{bingham,bojanic,galambos}. 

We say that a sequence of positive numbers $\lef \beta_n\ri$ is regularly
varying with exponent $\mu$ if the following holds: 
\[\beta_n\sim\gamma_n n^{\mu} \ \ \ \lef n\tv\ri\]
where $\lef\gamma_n\ri$ is slowly varying.

$\lef\beta_n\ri$ is regularly varying with exponent $\mu$ if and only if
 $\beta_{\lsz sn\rsz}/\beta_n\rightarrow s^{\mu}$ as $n\tv$ for all
 $s>0$, see Bingham \cite{bingham}.

\begin{prop}\label{prop4}
 Let $\lef \al_n\ri$, $\lef \beta_n\ri$ be nonnegative sequences such that
 $\lef \al_n \ri$ is slowly varying as $n\tv$, and $n^{-\lambda}
 \beta_n\rightarrow 1$ as $n\tv$ for some $\lambda>-1$. Then the
 following holds. 
\[\sum_{i=1}^n \al_i\beta_i\sim \al_n\sum_{i=1}^n \beta_i \quad\lef n\tv\ri.\]
\end{prop}

This is a consequence of the results of Bojani{\'c} and Seneta
\cite{bingham,bojanic}. 

\begin{prop}\label{prop5-6}
Let $\lef \al_n\ri$, $\lef \beta_n\ri$ be nonnegative sequences such that
$\lef \al_n \ri$ is regularly varying with exponent $\delta$.
\begin{itemize}
\item[\aaa] Suppose $\sum_{i=1}^n{\beta_i}=B_n$ is regularly varying with
exponent $\mu>0$, and $\mu+\delta>0$. Then
\[\sum_{i=1}^n \al_i\beta_i\sim \frac{\mu}{\delta+\mu}\,
\al_nB_n\quad\lef n\tv\ri.\]
\item[\bb] Suppose $\sum_{i=1}^n{\beta_i}=o(B_n)$, where $(B_n)$ is
regularly varying with exponent $\mu>0$, and $\mu+\delta>0$. Then 
\[\sum_{i=1}^n \al_i\beta_i=o\lef\al_nB_n\ri\quad\lef n\tv\ri.\] 
\end{itemize}
\end{prop}

\biz {\aaa} Suppose first that $\delta=0$, that is, $\lef\al_n\ri$ is slowly varying. By Bojani{\'c} and Seneta
\cite{bojanic}, for a nonnegative slowly varying sequence $\lef \al_n\ri$
there always exists another nonnegative sequence $\lef \al'_n\ri$ such
that $\al_n\sim  \al'_n$ as $n\tv$, and  
\begin{equation}\label{p1}
\lim_{n\tv} n\lef 1-\frac{\al'_{n-1}}{\al'_n}\ri=0.
\end{equation}
This implies that $\al_{n+1}/\al_n\rightarrow 1$ as $n\tv$. 

All sequences are nonnegative, hence we have
\begin{multline}\label{p1.5}
\sum_{i=1}^n \al'_i\beta_i=\al'_n\sum_{j=1}^n
\beta_j+\sum_{i=1}^{n-1}\lef \al'_i-\al'_{i+1}\ri\sum_{j=1}^i \beta_j\\ 
=\al'_n\sum_{j=1}^n \beta_j-\sum_{i=1}^{n-1}i\lef
1-\frac{\al'_{i}}{\al'_{i+1}}\ri\lef \frac{\al'_{i+1}}{i}\sum_{j=1}^i
\beta_j\ri
\end{multline} 
as $n\tv$.  

Sequence $\lef \al'_n\ri$ is slowly varying. By supposition,
$\sum_{j=1}^n\beta_j=\gamma_n n^{\mu}$ as $n\tv$ with some slowly
varying sequence $\lef\gamma_n\ri$, hence $i^{-1}\sum_{j=1}^i \beta_j=
\gamma_i\,i^{\mu-1}$ as $i\tv$. Since $\lambda=\mu-1>0$, by applying
Proposition \ref{prop4} we obtain that    
\[\sum_{i=1}^n \frac{\al'_{i+1}}{i}\sum_{j=1}^i \beta_j \sim
\sum_{i=1}^n \al'_{i+1} \gamma_i\, i^{\mu-1}\sim 
\al'_{n+1}\gamma_n\sum_{i=1}^n i^{\mu-1}\sim \frac{1}{\mu}\,\al_{n}
\gamma_n n^{\mu}\] 
as $n\tv$. Combining this with \eqref{p1} we get that the second
term on the right-hand side of \eqref{p1.5} is $o\lef \al_n \gamma_n
n^{\mu}\ri$ as $n\tv$.  

The first term is asymptotically equal to $\al_n \gamma_n n^{\mu}$ as
$n\tv$. Thus we get that  
\[\sum_{i=1}^n \al_i\beta_i\sim \sum_{i=1}^n \al'_i\beta_i\sim 
\al_n \gamma_n n^{\mu}\quad\lef n\tv\ri.\]

Next, let $\delta$ differ from $0$. Let $\alpha_n=\kappa_n
n^{\delta}$, and $B_n=\gamma_n n^{\mu}$ with slowly varying sequences
$(\kappa_n)$ and $(\gamma_n)$. We have
\begin{multline*}
\sum_{i=1}^n i^{\delta}\beta_i
=\sum_{i=1}^n i^{\delta}\lef B_i-B_{i-1}\ri
=n^{\delta}B_n+\sum_{i=1}^{n-1} \lef i^{\delta}-\lef i+1\ri^{\delta}\ri
B_i\\
=n^{\delta}B_n-\delta\sum_{i=1}^{n-1} i^{\delta-1} \lef 1+o\lef 1\ri\ri
B_i\\ 
=\gamma_nn^{\delta+\mu}- (1+o(1))\,\delta \sum_{i=1}^{n-1} 
\gamma_i i^{\delta+\mu-1}\quad\lef n\tv\ri.
\end{multline*} 
$\lef \gamma_n\ri$ is slowly varying, and $\lambda=\delta+\mu-1>-1$, thus
Proposition \ref{prop4} applies, and we obtain that 
\begin{multline*}
\sum_{i=1}^n i^{\delta} \beta_i \sim \gamma_n
n^{\delta+\mu}-\delta\gamma_n \sum_{i=1}^{n-1}  i^{\delta+\mu-1}\\ 
\sim \gamma_n n^{\delta+\mu}-\delta\gamma_n\,
\frac{n^{\delta+\mu}}{\delta+\mu}
= \frac{\mu}{\delta+\mu}\,\gamma_n
n^{\delta+\mu} \quad\lef n\tv\ri.               
\end{multline*} 
Let us apply the already proved particular case to $\lef\kappa_n\ri$ and
$\lef n^{\delta}\beta_n\ri$. Then we get that 
\[
\sum_{i=1}^n\al_i\beta_i=\sum_{i=1}^n \kappa_i i^{\delta} \beta_i
\sim \kappa_n\sum _{i=1}^n i^{\delta} \beta_i
\sim \frac{\mu}{\delta+\mu}\,\kappa_n \gamma_n n^{\delta+\mu}
=\frac{\mu}{\delta+\mu}\,\al_nB_n.\]
{\bb}\ \ We can
suppose that $B_n$ is increasing, since for every regularly varying
sequence with positive exponent one can find another, increasing one,
which is equivalent to it. Introduce
$\beta_n'=\beta_n+B_n-B_{n-1}$, with $B_0=0$. Then $\beta_n'\ge 0$, and 
\[\sum_{i=1}^n\beta_i'=\sum_{i=1}^n\beta_i+B_n\sim B_n,\]
hence it is regularly varying with exponent $\mu$. By part {\aaa} we have
\[\sum_{i=1}^n\al_i\beta_i'\sim\frac{\mu}{\delta+\mu}\,\al_nB_n,\]
and also
\[\sum_{i=1}^n\al_i\lef
B_i-B_{i-1}\ri\sim\frac{\mu}{\delta+\mu}\,\al_nB_n.\]
After subtraction we obtain that
$\sum_{i=1}^n\al_i\beta_i=o(\al_nB_n)$. \qed

\begin{prop}\label{l5}
Let $a_1,a_2\ldots$ and $b_1,b_2,\ldots$ be nonnegative numbers
satisfying 
\[\lim_{n\tv}\frac 1n \sum_{i=1}^n a_i=K<\infty,\quad \lim_{n\tv}nb_n=1.
\]
Then
\begin{itemize}
\item[\aaa] $\exp\lef \sum_{i=1}^n a_ib_i\ri$ is regularly varying with
  exponent $K$ as $n\tv$; 
\item[\bb] $\exp\lef \sum_{i=1}^n a_i^2b_i^2s_i\ri$ is slowly varying as
  $n\tv$ for every bounded sequence of real numbers $\lef s_n\ri$.
\end{itemize}
\end{prop}
\biz {\aaa} Let $m=[tn], \ t>1$. We have
\begin{multline*}
\sum_{i=n+1}^m a_i b_i=\sum_{i=n+1}^m\frac{a_iq_i}{i}+
\sum_{i=n+1}^m\frac 1i\lsz i(K+r_i)-(i-1)(K+r_{i-1})\rsz\\ 
=\sum_{i=n+1}^m\frac{a_iq_i}{i}+K\sum_{i=n+1}^m\frac 1i
+\sum_{i=n+1}^m\lef r_i-r_{i-1}\ri+\sum_{i=n+1}^m \frac{r_{i-1}}{i}. 
\end{multline*}
The first sum on the right-hand side tends to $0$, since 
\begin{equation}\label{becsles1}
\left\vert\sum_{i=n+1}^m \frac{a_iq_i}{i}\right\vert\leq
\sum_{i=n+1}^m\frac{a_i\left\vert q_i\right\vert}{i}\leq 
\frac tm\sum_{i=n+1}^m a_i\left\vert q_i\right\vert=o\lef 1\ri.
\end{equation}
The second sum is $K\log t+o\lef 1\ri$, the third one is $r_m-r_{n-1}=o\lef
1\ri,$ and the last one also converges to $0$.
 
\noindent{\bb} Let $m=[tn],\ t>1$. Now we have 
\[\sum_{i=n+1}^m a_i^2b_i^2 s_i= 
\sum_{i=n+1}^m \frac{a_i}{i}\cdot\frac{\lef 1+q_i\ri^2a_is_i}{i}.\]
By supposition
\[\frac{a_i}{i}=\frac{K+ir_i-\lef i-1\ri
  r_{i-1}}{i}=\frac{K}{i}+r_i-r_{i-1}+\frac{r_{i-1}}{i}\rightarrow 0
\quad\lef i\tv\ri, \] 
hence
\[q_i'=\frac{\lef 1+q_i\ri^2a_is_i}{i}\rightarrow 0\ \ \lef i\tv\ri.\] 
We can complete the proof similarly to  \eqref{becsles1}.
\[\left\vert\sum_{i=n+1}^m a_i^2b_i^2 s_i\right\vert= 
\left\vert\sum_{i=n+1}^m \frac{a_iq_i'}{i}\right\vert
\leq\frac tm\sum_{i=n+1}^m a_i\left\vert q_i'\right\vert=o\lef 1\ri.
\qed\]

\section{Proofs}\label{Proofs}

For sake of convenience, instead of $X\cs \lsz n, d\rsz$, we consider
the number $Z\cs[n,d]$ of selected vertices with degree greater than
or equal to $d$. That is, for $n\geq 1$ and $d\geq m$ let 
\begin{equation}\label{e-2}
Z\cs\lsz n, d\rsz= \left\vert \left\lbrace  v\in S_n:\ \deg_n\left(
v\right)\geq d\right\rbrace \right\vert
=\sum_{j=d}^n X\cs\lsz n, j\rsz.
\end{equation}

We also need the following notations.
\[J\cs\lsz n, d\rsz=\sum_{j=d}^n I\cs\lsz n, j\rsz,\quad
J\lsz n, d\rsz= \sum_{j=d}^n I\lsz n, j\rsz.\] 

First we show that Theorem \ref{t1} is implied by the following
proposition. For all $d\geq m$ we have $Z\cs\lsz n, d\rsz\sim
z_d\left\vert S_n\right\vert$ a.s. as $n\tv$ with some positive
constants $z_d$. In addition,   
\begin{equation}\label{rekurzio}
z_m=1,\quad z_d=\frac{z_{d-1}\frac{k_{d-1}}{c_{d-1}}+\alpha
  \sum_{j=d}^{\infty}q_j}{\alpha+\frac{k_{d-1}}{c_{d-1}}}
\quad (d\geq m+1).
\end{equation}

It is clear that 
\[X\cs\lsz n, d\rsz=Z\cs\lsz n, d\rsz-Z\cs\lsz n, d+1\rsz \quad \lef
n\geq 1,\ d\geq m\ri, \] 
hence
\[X\cs\lsz n, d\rsz = \lef z_d-z_{d+1}\ri\left\vert S_n\right\vert +o\lef
\left\vert S_n\right\vert\ri\]
a.s. as $n\tv$. Thus the limits 
\[\lim_{n\rightarrow\infty}\frac{X\cs\left[ n, d\right]}{\left\vert
    S_n\right\vert}=x_d\] 
exist for all $d\geq m$ almost surely, and $x_d=z_d-z_{d+1}$ for all
$d\geq m$. 

It is easy to derive the recursive equations for $x_d=z_d-z_{d+1}$
from $z_m=1$ and equation \eqref{rekurzio}. The denominators are
positive, because Conditions \ref{c1}, \ref{c6}, and \ref{d3}
guarantee that $c_d$ is nonnegative, $\al$ is positive, and $k_d$ is
positive.   

It is also easy to check that sequence $\lef x_d\ri$ is a probability
distribution. We have 
\begin{equation*}
 x_m \al+x_m \frac{k_m}{c_m}= \al q_m,
\end{equation*}
and
\begin{equation*}
 x_d \al +x_d\, \frac{k_d}{c_d}=x_{d-1}\,\frac{k_{d-1}}{c_{d-1}}+\al
 q_d  \quad\lef d\geq m+1\ri.
\end{equation*}

Summing up the equations above we get that
\[\sum_{d=m}^{\infty} x_d=\sum_{d=m}^{\infty} q_d=1,\]
since, by Conditions \ref{d3} and \ref{d4}, $\al>0$ and the sequence
$(q_d)$ is a probability distribution. 

The next step is solving the recursion for $\lef x_d\ri$. 
Set
\[t_d=\frac{k_d}{c_d},\quad a_d=\prod_{i=m}^{d-1}
\frac{t_i+\alpha}{t_i}\quad \lef d\geq m\ri.\]  
It is easy to check that the recursive equations of Theorem
\ref{t1} are satisfied by the sequence 
\begin{equation*}
 x_d=\frac{1}{t_d+\alpha}\;\sum_{i=0}^d q_i\alpha\,\prod_{j=i}^{d-1}
 \frac{t_j}{t_j+\alpha}\quad \lef d\geq m\ri. 
\end{equation*}
By Condition \ref{c2}, $c_d\sim K\cdot d^{-\gamma}$ holds as $d\tv$,
and by Condition \ref{c4} the sequence $\lef p_j\ri$ is exponentially
decreasing. Hence it follows, as $d\tv$, that 
\begin{align*}
k_d&=-\sum_{j=m}^d \lef c_j-p_j\ri\sim -K\cdot
\frac{d^{-\gamma+1}}{-\gamma+1};\\
t_d&=\frac{k_d}{c_d}\sim\frac{-K\cdot\tfrac{d^{-\gamma+1}}{-\gamma+1}}
{K\cdot d^{-\gamma}}=\frac{d}{\gamma-1};\\
a_d&=\prod_{i=0}^{d-1} \lef 1+\frac{\al}{t_i}\ri\sim
\prod_{i=0}^{d-1} \lef 1+\frac{\al\lef \gamma-1\ri}{i}\ri\sim 
K'\cdot d^{\alpha\lef \gamma-1\ri}
\end{align*}
for some $K'>0$. 
By Condition \ref{d4} the sequence $\lef q_d\ri$ is exponentially
decreasing, thus the series in the expression
\begin{equation*}
 x_d=\frac{1}{a_d\lef t_d+\al\ri}\sum_{i=m}^da_iq_i\alpha
\end{equation*}
converges. Using the asymptotics of $\lef a_d\ri$ and $\lef t_d\ri$ we get
that  
\begin{equation*}
x_d=\frac{1}{a_d\lef t_d+\al\ri}\sum_{i=0}^da_iq_i\alpha\sim L\cdot
d^{-\alpha\lef \gamma-1\ri-1} 
\end{equation*}
for some $L>0$.

Consequently, the degree distribution constrained on the set of
selected vertices decays polynomially, and the new characteristic
exponent is determined by $\al$ and $\gamma$, namely,
$\gamma\cs=\alpha\lef \gamma-1\ri+1$, as stated.

Therefore Theorem \ref{t1} is indeed a consequence of \eqref{rekurzio}.

Let us continue with the proof of \eqref{rekurzio}. We proceed by
induction on $d$.  

The case $d=m$ is obvious, because the initial degree is never less
than $m$, and the degree of a vertex cannot decrease, thus every
vertex in $S_n\setminus S_0$ has at least $m$ edges.

Suppose that 
\begin{equation*}
Z\cs\lsz n, d-1\rsz\sim z_{d-1} \left\vert
S_n\right\vert\quad\lef n\tv\ri
\end{equation*}  
holds for some $z_{d-1}>$ and $d\geq m+1$ almost surely. 

First we determine the expected number of vertices of degree $\geq d$
in $S_{n+1}$, given $\F_n$, for $n\geq 1$. Every vertex
in $S_n$ counts if its degree is at least $d$ in $G_n$, or if its
degree is equal to $d-1$ in $G_n$ and it gets a new edge from
$v_{n+1}$. The new vertex $v_{n+1}$ counts if it falls into $S_{n+1}$
and its degree is $\geq d$ in $G_{n+1}$. Thus  the
following equality holds for every $n\geq 1$.  
\begin{equation}\label{e-1}
 Z\cs\lsz n+1, d\rsz=Z\cs\lsz n, d\rsz + Y\cs\lsz n, d-1\rsz+J\cs\lsz
 n+1, d\rsz. 
\end{equation}
Taking conditional expectations with respect to $\F_n$ we obtain that
\begin{multline}\label{e-0.5}
E\lef \left.Z\cs\lsz n+1, d\rsz\right\vert\F_n\ri\\
=Z\cs\lsz n, d\rsz+E\lef\left.Y\cs\lsz n, d-1\rsz \right\vert\F_n\ri
+E\lef\left.J\cs\lsz n+1, d\rsz\right\vert\F_n\ri. 
\end{multline}
By Condition \ref{c3}, vertices of the same degree are connected to
$v_{n+1}$ with the same conditional probability. This 
implies that  
\begin{equation}\label{e0}
 E\lef \left.\frac{Y\cs\lsz n, d\rsz}{X\cs\lsz n, d
 \rsz}\right\vert\F_n\ri=E\lef\left. \frac{Y\lsz n, d\rsz}{X\lsz n, d
 \rsz}\right\vert\F_n\ri,\quad(n\geq 1).
\end{equation} 
$X\lsz n, d\rsz$ may be equal to zero, then $Y\lsz n, d\rsz=0$
as well. We will consider all quotients of the form $0/0$ as $1$. 

The middle term on the right-hand side of \eqref{e-0.5} can be
transformed by the help of \eqref{e0}. 
\begin{multline}\label{e1}
E\lef \left.Z\cs\lsz n+1,d\rsz\right\vert\F_n\ri\\
=Z\cs\lsz n,d\rsz+X\cs\lsz n, d-1\rsz\frac{E\lef\left.Y\lsz n,d-1\rsz
\right\vert\F_n\ri}{X\lsz n, d-1\rsz}+E\lef\left.J\cs\lsz
n+1,d\rsz\right\vert\F_n\ri. 
\end{multline}
By \eqref{e-2}, $X\cs \lsz n, d-1\rsz=Z\cs\lsz n,
d-1\rsz-Z\cs\lsz n, d\rsz$, hence from equation \eqref{e1} we
obtain that   
\begin{multline}\label{e1.1}
E\lef \left.Z\cs\lsz n+1,d\rsz\right\vert\F_n\ri=Z\cs\lsz n, d\rsz\lef
1-\frac{E\lef\left.Y\lsz n, d-1\rsz\right\vert\F_n\ri}{X\lsz n,
  d-1\rsz}\ri\\
+Z\cs\lsz n, d-1\rsz\frac{E\lef\left.Y\lsz n,
d-1\rsz\right\vert\F_n\ri}{X\lsz n, d-1\rsz}+E\lef\left.J\cs\lsz n+1,
d\rsz\right\vert\F_n\ri 
\end{multline}
for all $n\geq 1$.

For $i\geq 1$ define 
\begin{equation*}
 b\lsz i, d\rsz=\begin{cases}
                 1 & \text{if } X\lsz i, d\rsz=0;\\
\lef 1-\frac{E\lef\left.Y\lsz i,
 d\rsz\right\vert\F_i\ri}{X\lsz i, d\rsz}\ri^{-1} & \text{ if } X\lsz i, d\rsz>0.
                \end{cases}
\end{equation*}

Set $c\lsz 1, d\rsz=1$ and for $n\geq 2$ define
\begin{equation}\label{e2}
 c\lsz n, d\rsz=\prod_{i=1}^{n-1} b\lsz i, d\rsz. 
\end{equation}
Then for $n$ large enough we have 
\begin{equation*}
 \frac{c\lsz n, d\rsz}{c\lsz n+1, d\rsz}=\lef 1-\frac{E\lef\left.Y\lsz n,
 d\rsz\right\vert\F_n\ri}{X\lsz n, d\rsz}\ri. 
\end{equation*}

For several particular models it is quite easy to compute the
conditional expectations $E\lef\left.Y\lsz i, d-1\rsz\right\vert\F_i\ri$,
and hence, to determine the asymptotics of $c\lsz n, d\rsz$. In the
present general case the conditional expectation is not
specified. However, as the following sequence of lemmas shows, the
asymptotics of the partial sums can be described, and one can  
calculate the asymptotics of $c\lsz n, d\rsz$. The proof of the lemmas
will be postponed to the second part of this section. We emphasize
that in the lemmas the induction hypothesis is assumed all along.

Consider the partial sums
\begin{equation*}\label{a1}
S\lsz n, d\rsz = \sum_{i=1}^n  E\lef\left. Y\lsz i,
  d\rsz\right\vert\F_i\ri\quad \lef n\geq 1\ri. 
\end{equation*}
\begin{lemma} \label{l1}
For all $d\geq m$ we have 
\begin{equation}
 S\lsz n, d\rsz = \sum_{i=1}^n E\lef\left. Y\lsz i,
 d\rsz\right\vert\F_i\ri= k_d\cdot n +o\lef n\ri \ \ \lef n\tv\ri  
\end{equation}
with probability $1$.
\end{lemma}

\begin{remark}It is clear from the definition that $S\lsz n, d\rsz$ is
  nonnegative, hence Lemma \ref{l1} immediately implies $k_d\geq 0$ for
  all $d\geq m$ (cf.\ Condition \ref{c6}).
\end{remark}  

\begin{lemma}\label{l2}
\begin{equation*} 
 c\lsz n, d\rsz\sim a\lsz n, d\rsz\cdot n^{k_d/c_d} \quad\lef n\tv \ri,
\end{equation*}
a.s. for all $d\geq m$, where $a\lsz n, d\rsz$ is positive and slowly
varying as $n\tv$.  
\end{lemma}

By equation \eqref{e1.1}, the process
\begin{equation}\label{vnd}
V\lsz n, d\rsz= c\lsz n, d-1\rsz Z\cs\lsz n, d\rsz\quad (n\geq 1)
\end{equation}
is a submartingale. Let $A\lsz n,
d\rsz$ denote the increasing process in the Doob decomposition of
$V\lsz n,d\rsz$; it is given by
\begin{multline}\label{e3}
 A\lsz n, d\rsz=\sum_{i=1}^n c\lsz i+1, d-1\rsz  Z\cs\lsz i,
 d-1\rsz\frac{E\lef\left.Y\lsz i, d-1\rsz\right\vert\F_i\ri}{X\lsz i,
 d-1\rsz}\\+\sum_{i=1}^n c\lsz i+1, d-1\rsz E\lef\left.J\cs\lsz i+1, d\rsz\right\vert\F_i\ri. 
\end{multline}
First we describe the asymptotics of $A\lsz n, d\rsz$.
\begin{lemma}\label{l3} 
Suppose that $Z\cs\lsz n, d-1\rsz\sim   z_{d-1}\left\vert
  S_n\right\vert$ holds a.s. for some $d\geq m+1$, as $n\tv$, then    
\begin{equation*}
 A\lsz n, d\rsz\sim \frac{z_{d-1}\frac{k_{d-1}}{c_{d-1}}+\alpha
 \sum_{j=d}^{\infty}q_j}{\alpha+\frac{k_{d-1}}{c_{d-1}}}\,a\lsz n, d\rsz
 \zeta_nn^{\alpha+k_{d-1}/c_{d-1}}\quad\textrm {a.s.} 
\end{equation*}
\end{lemma}
Next, we compute an upper bound for the conditional variances. Define 
\begin{equation*}
 B\lsz n, d\rsz =\sum_{i=2}^n \textrm{Var}\lef\left. V\lsz i,
 d\rsz\right\vert\F_{i-1}\ri \quad\lef n\geq 2\ri. 
\end{equation*}
\begin{lemma}\label{l4}
Suppose that $Z\cs\lsz n, d-1\rsz\sim   z_{d-1}\left\vert
S_n\right\vert$ holds a.s. for some $d\geq m+1$, as $n\tv$, then
$B\lsz n, d\rsz^{1/2}\log B\lsz n, d\rsz=O\lef A\lsz n, d\rsz\ri$. 
\end{lemma}

Therefore Proposition \ref{prop1} implies that $V\lsz n, d\rsz\sim
A\lsz n, d\rsz$ almost surely as $n\tv.$ Finally, by Lemma \ref{l2}
and Lemma  \ref{l3} we obtain the asymptotics
\begin{equation*}
 Z\cs\lsz n, d\rsz \sim \frac{z_{d-1}\frac{k_{d-1}}{c_{d-1}}+\alpha
 \sum_{j=d}^{\infty}q_j}{\alpha+\frac{k_{d-1}}{c_{d-1}}}\,\zeta_n
 n^{\alpha} \ \ \ \lef n\tv\ri. 
\end{equation*}
Consequently, we have
\begin{equation*}Z\cs\lsz n, d\rsz\sim z_d\zeta_nn^{\al} \quad\lef
  n\tv\ri
\end{equation*}
with
\begin{equation}
z_d=\frac{z_{d-1}\frac{k_{d-1}}{c_{d-1}}+\alpha
  \sum_{j=d}^{\infty}q_j}{\alpha+\frac{k_{d-1}}{c_{d-1}}}\,. 
\end{equation}
 
The size of $S_n$ is asymptotically equal to $\zeta_nn^{\al}$ by
Condition \ref{d3}. Thus the proof of \eqref{rekurzio} can
be completed by using Lemmas \ref{l1}--\ref{l4}.\qed

\medskip  
Now we continue with the proofs of Lemmas \ref{l1}--\ref{l4}.

\medskip

\textbf{Proof of Lemma \ref{l1}.} 
Similarly to equation \eqref{e-1}, but considering all vertices,
we see that 
\begin{equation*}\label{eq1mas}
 X\lsz i+1, j\rsz=X\lsz i, j\rsz -Y\lsz i, j\rsz+Y\lsz i,
 j-1\rsz+I\lsz i+1, j\rsz
\end{equation*}
for every $i\geq 0$ and $j\geq m$. 
Adding up for $i=1,\ldots,n$  we
obtain that 
\begin{equation}\label{e11}
 X\lsz n+1, j\rsz-X\lsz 1, j\rsz=-\sum_{i=1}^n Y\lsz i, j\rsz+\sum_{i=1}^n Y\lsz i,
 j-1\rsz+\sum_{i=2}^{n+1} I\lsz i, j\rsz 
\end{equation}
for every $j\geq m $ and $n\geq 1$.
By Conditions \ref{c1} and \ref{c4}, from \eqref{e11} it follows that
\begin{equation*}
\sum_{i=1}^n Y\lsz i, j-1\rsz-\sum_{i=1}^n Y\lsz i, j\rsz=\lef
 c_j-p_j\ri\cdot n+o\lef n\ri
\end{equation*}
holds almost surely, as $n\tv$, for every $j\geq m$.
Adding this up for $j=m,\ldots, d$ we get
\begin{equation}\label{e12}
 \sum_{i=1}^n Y\lsz i, d\rsz=-\sum_{j=m}^d\lef c_j-p_j\ri\cdot n+o\lef
 n\ri=k_d\cdot n+o\lef n\ri \end{equation} 
a.s., as $n\tv$. Therefore it is sufficient to prove that 
\begin{equation*}\frac{1}{n}\sum_{i=1}^n \lef Y\lsz i, d\rsz -E\lef
\left. Y\lsz i, d\rsz\right\vert\F_i\ri\ri\rightarrow 0\quad \lef n\tv\ri.
\end{equation*} 
Fix $d\geq m$, and for $n\geq 1$ let 
$M_n=\sum_{i=1}^n \lef Y\lsz i, d\rsz -E\lef\left. Y\lsz i,
  d\rsz\right\vert\F_i\ri\ri$, $\mathcal G_n=\F_{n+1}$. It is clear that
$\lef M_n, \mathcal G_n\ri$ is a martingale. Using Condition \ref{c5} we
will derive an upper bound for the corresponding increasing process
$A_n$ introduced in Proposition \ref{prop0}. 
\begin{multline}\label{l12}
A_n=\sum_{i=1}^n\textrm{Var}\lef\left. Y\lsz i,d\rsz\right\vert\F_i\ri\leq
\sum_{i=1}^n E\lef \left.Y\lsz i, d\rsz^2\right\vert\F_i\ri\leq\\ 
\leq\sum_{i=1}^n C_i E\lef\left. Y[i,d]\right\vert\F_i\ri+
\sum_{i=1}^n E\lef \left.Y[i,d]^2I\lef Y[i,d]>C_i\ri\right\vert\F_i\ri\\
\leq\sum_{i=1}^n C_iE\lef\left. Y\lsz i,d\rsz\right\vert\F_i\ri+
\sum_{i=1}^n E\lef Z_d^2I\lef Z_d>C_i\ri\ri
\end{multline}
for any $C_i>0$. Fix $\e>0$ such that $\kappa=E\lef e^{\e Z_d}\ri$ is
finite, and for $i\geq 3$ choose $C_i=\frac 2{\e}\log i$. The function
$z\mapsto z^2e^{-\e z}$ is decreasing for $z>\frac{2}{\e}$, hence
$z^2e^{-\e z}\leq C_i^2e^{-\e C_i}$ for $z> C_i$. This implies 
\begin{multline}\label{l13}
E\lef Z_d^2I\lef Z_d>C_i\ri\ri\leq C_i^2e^{-\e C_i}E\lef e^{\e Z_d}I\lef Z_d>
C_i\ri\ri\\
\leq (2/\e)^2\lef \log i\ri^2i^{-2}\kappa.
\end{multline}
The infinite sum of these terms converges, thus the second sum on the
right-hand side of \eqref{l12} is bounded for fixed $d$.  

On the other hand,  $Y\lsz i, d\rsz \leq \left\vert V_i\right\vert
\leq i+l$ follows from the definition, therefore
\[\sum_{i=3}^n C_iE\lef \left. Y\lsz i, d\rsz \right\vert\F_i\ri\leq
\sum_{i=3}^n \frac 2{\e}\log i \cdot \lef i+l\ri=O\lef n^2\log n\ri.\] 
Thus $A_n=O\lef n^2\log n\ri$. This bound can be further improved as
follows. Applying Proposition \ref{prop0} to the martingale $\lef
M_n\ri$ we get that $M_n=O\lef n^{1+\eta}\ri$ a.s. for all $\eta>0$.
Equation \eqref{e12} implies that $\sum_{i=1}^n Y\lsz i, d\rsz=O\lef
 n\ri$, therefore
\begin{multline*}
 \sum_{i=1}^n C_iE\lef \left. Y\lsz i,d\rsz \right\vert\F_i\ri\leq
 C_n\sum_{i=1}^n E\lef \left. Y\lsz i,d\rsz \right\vert\F_i\ri\\
 =C_n\lef\sum_{i=1}^n Y\lsz i,d\rsz -M_n\ri=
 O\lef n^{1+\eta}\log n\ri.
\end{multline*}
We obtain that $A_n=O\lef n^{1+\eta}\log n\ri$. Hence by Proposition
\ref{prop0} we have $M_n=o\lef n^{\frac 12+\eta}\log n\ri$ a.e. on the
event $\lk A_{\infty}=\infty\rk$, for all $\eta>0$. Therefore
$M_n=o(n)$ holds almost surely, and this completes the proof of 
Lemma \ref{l1}.\qed  

\medskip
\textbf{Proof of Lemma \ref{l2}.}
Fix an arbitrary $d\geq m$. Lemma \ref{l1} and the induction
hypothesis imply that 
\[\frac{E\lef\left.Y\lsz n, d\rsz\right\vert\F_n\ri}{X\lsz n,
  d\rsz}=\frac{S\lsz n, d\rsz-S\lsz n-1, d\rsz}{X\lsz n,
  d\rsz}\rightarrow 0\quad\lef n\tv\ri.\] 
Thus in \eqref{e2} we can apply the approximation
$1-x=e^{-x+O\lef x^2\ri}$ $(x\to 0)$. Set
\[a_i=\frac{S\lsz i, d\rsz-S\lsz i-1, d\rsz}{c_d}\,, \quad 
b_i=\frac{c_d}{X\lsz i, d\rsz}\]
if $X\lsz i, d\rsz>0$, and $a_i=b_i=0$ otherwise. Then $a_i$ and $b_i$ are
nonnegative. In addition, 
\[\frac 1n\sum_{i=1}^n a_i=\frac 1{c_dn}\sum_{i=1}^n\lef S\lsz i,d\rsz-
S\lsz i-1, d\rsz\ri=\frac{S\lsz n,d\rsz}{c_dn}\to\frac{k_d}{c_d}\]  
as $n\tv$, by Lemma \ref{l1}. According to the induction hypothesis,
$X\lsz n, d\rsz\sim c_d\cdot n$, which implies that
$n b_n\rightarrow 1$ as $n\tv$.
Therefore Proposition \ref{l5} applies to the sequences $\lef a_n\ri$ and
$\lef b_n\ri$ with $K=k_d/c_d$. Thus, due to part {\aaa},  
\[\exp\lef \sum_{i=1}^n a_ib_i\ri=\exp\lef\sum_{i=1}^n \frac{S\lsz i,
  d\rsz-S\lsz i-1, d\rsz}{X\lsz i, d\rsz} \ri \]  
is regularly varying with exponent $K$. 

The remainder terms produce a slowly varying function, because by part
{\bb} of Proposition \ref{l5} we get that
\[\exp\lef \sum_{i=1}^n a_i^2b_i^2s_i\ri=
\exp\lef\sum_{i=1}^n \lef\frac{S\lsz i,d\rsz-S\lsz i-1, d\rsz}{X\lsz
i,d\rsz}\ri^{\!2}s_i\ri \] 
is slowly varying supposed the sequence $\lef s_i\ri$ is bounded.

From these the asymptotics of $c\lsz n, d\rsz$ readily follows.\qed

\medskip

\textbf{Proof of Lemma \ref{l3}.} By \eqref{e3}, $A[n,d]=A_1+A_2$,
where
\begin{align*}
A_1&=\sum_{i=1}^n c\lsz i+1, d-1\rsz Z\cs\lsz i, d-1\rsz 
\frac{E\lef\left.Y\lsz i,d-1\rsz\right\vert\F_i\ri}{X\lsz i,d-1\rsz}\,,\\
A_2&=\sum_{i=1}^n c\lsz i+1, d-1\rsz E\lef\left.J\cs\lsz i+1,d\rsz
\right\vert\F_i\ri.
\end{align*}

Here we already know the asymptotics of $S\lsz n, d-1\rsz$ and $c\lsz
n, d-1\rsz$ from Lemmas \ref{l1} and \ref{l2}. In addition, $Z\cs\lsz n, 
d-1\rsz\sim z_{d-1}\left\vert S_n\right\vert\sim
z_{d-1}\zeta_nn^{\al}$ a.s., due to the induction
hypothesis and Condition \ref{d3}. It is clear from the
definition that $Y\lsz n, d-1\rsz$ is nonnegative, thus we have 
\begin{multline*}
A_1\sim\sum_{i=1}^n a\lsz i+1,d-1\rsz \,i^{k_{d-1}/c_{d-1}}\,z_{d-1}
\,\zeta_i\,i^\al \,\frac{1}{c_{d-1}i}\,E\lef\left.Y\lsz i,d-1\rsz
\right\vert\F_i\ri\\
=\frac{z_{d-1}}{c_{d-1}}\sum_{i=1}^n a\lsz i+1, d-1\rsz\,\zeta_i\,
i^{k_{d-1}/c_{d-1}+\al-1}\,E\lef\left.Y\lsz i, d-1\rsz\right\vert\F_i\ri.
\end{multline*}

Let us apply part {\aaa} of Proposition \ref{prop5-6} in the following
setting.  
\[\al_n=a\lsz n+1, d-1\rsz\zeta_n,\quad
\beta_n=E\lef\left.Y\lsz n, d-1\rsz\right\vert\F_n\ri\]
for $n\geq 1$; 
\[\delta=\frac{k_{d-1}}{c_{d-1}}+\al-1.\]

Condition \ref{d3} and Lemma \ref{l2} guarantee that $\al_n$ is slowly
varying. Furthermore, $\sum_{i=1}^n \beta_i=S\lsz n, d-1\rsz\sim
k_{d-1} n$ as $n\tv$ with $k_{d-1}>0$, hence $\gamma_n=k_{d-1}$ and
$\mu=1$ satisfy the conditions. Finally,
$\mu+\delta=k_{d-1}/c_{d-1}+\al>0$, because 
$c_{d-1}$, $k_{d-1}$ and $\al$ are positive due to Conditions
\ref{c1}, \ref{c6}, and \ref{d3}.  

Applying Proposition \ref{prop5-6} we obtain that 
\[A_1\sim\frac{z_{d-1}k_{d-1}}{c_{d-1}}\cdot\frac{1}
{\al+\frac{k_{d-1}}{c_{d-1}}}\cdot a\lsz n,d-1\rsz\,\zeta_n\,
n^{k_{d-1}/c_{d-1}+\al}\] 
almost surely as $n\tv$, where $z_{d-1}, k_{d-1}$ and $c_{d-1}$ are
positive. 

Now we examine the second term in $A[n,d]$. Since 
\[J\cs[i+1,d]=I\cs(i+1)-\sum_{j=m}^{d-1}I\cs[i+1,j],\] 
we have
\[E\lef\left.J\cs\lsz i+1,d\rsz\right\vert\F_i\ri=
E\lef\left.I\cs(i+1)\right\vert\F_i\ri-
\sum_{j=m}^{d-1}E\lef\left.I\cs[i+1,j]\right\vert\F_i\ri.\]
Hence by Lemma \ref{l2}  
\begin{multline*}
A_2\sim\sum_{i=1}^n a\lsz i+1, d-1\rsz i^{k_{d-1}/c_{d-1}}\times\\
\times\lef E\lef\left. I\cs\lef i+1\ri\right\vert\F_i\ri-
\sum_{j=m}^{d-1} E\lef\left.I\cs\lsz i+1,j\rsz\right\vert\F_i\ri\ri. 
\end{multline*}

Set $\al_n=a\lsz n+1, d-1\rsz$, $\delta=k_{d-1}/c_{d-1}$, and  
$\beta_n=E\lef\left.I\cs\lef n+1\ri\right\vert\F_n\ri$. By Proposition
\ref{prop2} and  Condition \ref{d3} we have  
\begin{equation}\label{beta}
\sum_{i=1}^n\beta_i\sim\sum_{i=1}^n I\cs\lef i+1\ri=\left\vert
S_{n+1}\right\vert-I\cs\lef 1\ri\sim \zeta_n n^{\al}\quad \lef n\tv\ri.
\end{equation}
Thus we can apply part {\aaa} of Proposition \ref{prop5-6} with
$\mu=\al>0$. Assumption $\delta+\mu>0$ is satisfied. Therefore we get
that   
\[\sum_{i=1}^n c\lsz i+1, d-1\rsz I\cs\lef i+1\ri\sim
\frac{\al}{\al+\frac{k_{d-1}}{c_{d-1}}}\cdot a\lsz n+1,
d-1\rsz\,\zeta_n\,n^{\al+k_{d-1}/c_{d-1}} \]  
almost surely as $n\tv$.  

On the other hand, for a fixed $ j\leq d-1$ we have
\begin{multline*}
\sum_{i=1}^n c\lsz i+1, d-1\rsz E\lef\left.I\cs\lsz i+1,j\rsz
\right\vert\F_i\ri\\ 
\sim \sum_{i=1}^n a\lsz i+1, d-1\rsz
i^{k_{d-1}/c_{d-1}}E\lef\left.I\cs\lsz i+1, j\rsz\right\vert\F_i\ri 
\end{multline*}
by Lemma \ref{l2}. In this case $\al_n$ remains the same as before,
and we set $\beta_n=E\lef\left.I\cs\lsz n+1,j\rsz\right\vert\F_n\ri$. 
Using Condition \ref{d4} and equation \eqref{beta} we obtain that 
\begin{multline*}
\sum_{i=1}^n \beta_i =\sum_{i=1}^n E\lef\left.I\cs\lsz i+1,j\rsz\right
\vert\F_i\ri\\
=\lef q_j+o\lef 1\ri\ri\sum_{i=1}^n E\lef\left.I\cs\lef i+1\ri\right\vert\F_i\ri
=\lef q_j+o\lef 1\ri\ri\zeta_n n^{\al}
\end{multline*} 
almost surely as $n\tv$. 
Thus we can apply part {\aaa} or part {\bb} of Proposition \ref{prop5-6}
with $\mu=\al$, according that $q_j$ vanishes or it is positive. Then we
get that    
\begin{multline*}
\sum_{i=1}^n  c\lsz i+1, d-1\rsz E\lef\left.I\cs\lsz i+1,
   j\rsz\right\vert\F_i\ri\\=\frac{\al q_j+o\lef
   1\ri}{\al+\frac{k_{d-1}}{c_{d-1}}}\cdot a\lsz n,
 d-1\rsz\zeta_nn^{\al+k_{d-1}/c_{d-1}} 
\end{multline*} 
almost surely as $n\tv$. Hence we conclude that 
\begin{equation}\label{A2}
A_2\sim \lef 1-\sum_{j=m}^{d-1} q_j+o\lef 1\ri\ri
\frac{\al}{\al+\frac{k_{d-1}}{c_{d-1}}}\cdot a\lsz n,d-1\rsz\,\zeta_n
n^{\al+k_{d-1}/c_{d-1}}
\end{equation} 
almost surely as $n\tv$. 
Since  $\lef q_d\ri$ is a probability distribution by Condition \ref{d4},
it follows that $\lef 1-\sum_{j=m}^{d-1} q_j\ri=\sum_{j=d}^{\infty}q_j$. This 
completes the proof.\qed

\medskip
\textbf{Proof of Lemma \ref{l4}.}
From equation \eqref{vnd} it follows that 
\[B\lsz n, d\rsz=\sum_{i=2}^n\textrm{Var} \lef\left.V\lsz i,
  d\rsz\right\vert\F_{i-1}\ri=
\sum_{i=2}^n c\lsz i, d-1\rsz^2 \textrm{Var} \lef
 \left.Z\cs\lsz i, d\rsz\right\vert\F_{i-1}\ri.\] 
By equation \eqref{e-1} we have
\begin{multline*} 
 \textrm{Var}\lef \left.Z\cs\lsz i, d\rsz\right\vert\F_{i-1}\ri\leq
 E\lef\left.\lef Z\cs\lsz i, d\rsz-Z\cs\lsz i-1,
   d\rsz\ri^2\right\vert\F_{i-1}\ri\\  
=E\lef\left.\lef Y\cs\lsz i-1, d-1\rsz+J\cs\lsz i,
 d\rsz\ri^2\right\vert\F_{i-1}\ri \\ 
\leq 2 E\lef\left.Y\cs\lsz i-1,
  d-1\rsz^2\right\vert\F_{i-1}\ri+2E\lef\left.J\cs\lsz i,
  d\rsz\right\vert\F_{i-1}\ri. 
\end{multline*}
Hence
\begin{multline}\label{sum}
B\lsz n, d\rsz\leq 2\sum_{i=2}^n c[i,d-1]^2\,E\lef\left.Y\cs\lsz i-1,
d-1\rsz^2\right\vert\F_{i-1}\ri\\
\qquad\qquad+2\sum_{i=2}^n
c[i,d-1]^2\,E\lef\left.J\cs[i,d]\right\vert\F_{i-1}\ri 
=2B_1+2B_2.
\end{multline}

We will estimate $B_1$ and $B_2$ separately. 

Similarly to the proof of Lemma \ref{l1}, fix a positive $\e>0$ such that
$\kappa=E\lef e^{\e Z_{d-1}}\ri<\infty$, and set $C_i=\frac{2}{\e}\log
i$. Using Condition \ref{c5} and inequality $Y\cs\lsz i,
d\rsz\leq Y\lsz i,d\rsz$ one can see that     
\begin{multline*}
E\lef\left.Y\cs\lsz i-1,d-1\rsz^2\right\vert\F_{i-1}\ri\leq
C_i E\lef\left.Y\cs\lsz i-1, d-1\rsz\right\vert\F_{i-1}\ri\\
+E\lef Z_{d-1}^2I\lef Z_{d-1}>C_i\ri\ri
\end{multline*}
holds. For estimating the first term on the right-hand side we make
use of equation \eqref{e0}.  
\begin{multline*}
E\lef\left.Y\cs\lsz i-1,
  d-1\rsz\right\vert\F_{i-1}\ri\\=
\frac{E\lef\left.Y\lsz i-1, d-1
    \rsz\right\vert\F_{i-1}\ri}{X\lsz i-1, d-1\rsz}\,X\cs\lsz i-1,d-1
\rsz\\ 
\leq\frac{E\lef\left.Y\lsz i-1, d-1
    \rsz\right\vert\F_{i-1}\ri}{X\lsz i-1, d-1\rsz}\left\vert
  S_{i-1}\right\vert. 
\end{multline*}

To the second term we can apply \eqref{l13}; it is $O\lef(\log
i)^2i^{-2}\ri$.

From all these we obtain that
\begin{multline*}
B_1\leq\sum_{i=2}^n
c[i,d-1]^2C_i\,\frac{E\lef\left.Y[i-1,d-1]\right\vert\F_{i-1}\ri}
{X[i-1,d-1]}\left\vert S_{i-1}\right\vert\\
\qquad\qquad+O\lef\sum_{i=2}^n c[i,d-1]^2(\log i)^2i^{-2}\ri
\end{multline*} 
Note that the second sum is convergent here. In the first sum $c[i,d-1]$
  can be estimated by Lemma \ref{l2}, $\left\vert S_{i-1}\right\vert$
  by Condition \ref{d3}, and $X[i-1,d-1]$ by Condition \ref{c1}. In
  this way we obtain that  
\[c[i,d-1]^2\,C_i\,\frac{\left\vert S_{i-1}\right\vert}{X[i-1,d-1]}\]
is regularly varying with exponent $\delta=2k_{d-1}/c_{d-1}+\al-1$. On the
other hand, by Lemma \ref{l1} the sum of 
$E\lef\left.Y[i-1,d-1]\right\vert\F_{i-1}\ri$ is regularly varying with
exponent $1$. Therefore part {\aaa} of Proposition \ref{prop5-6} implies
that 
\[B_1=O\lef a[n,d-1]^2(\log n)^2\zeta_n\,n^{\al+2k_{d-1}/c_{d-1}}\ri.\]

For $B_2$ let us apply part {\aaa} of Proposition \ref{prop5-6} with
$\al_i=c[i,d-1]$ and 
$\beta_i=c[i,d-1] E\lef\left.J\cs[i,d]\right\vert\F_{i-1}\ri$. The regular
variation of $\sum\beta_i$ has already been proven in
\eqref{A2}. Thus,
\[B_2=O\lef a[n,d-1]^2\,\zeta_n\,n^{\al+2k_{d-1}/c_{d-1}}\ri.\]

Returning to \eqref{sum} we conclude that
\[B\lsz n, d\rsz =O\lef n^{\al+2 k_{d-1}/c_{d-1}+\eta}\ri\]
for all $\eta>0$.
Consequently, 
\[B\lsz n, d\rsz^{1/2} \log B\lsz n, d\rsz =O\lef
n^{\al/2+k_{d-1}/c_{d-1}+\eta}\ri\quad\lef n\tv\ri.\] 

Now the proof can be completed by comparing this with Lemma \ref{l3}.\qed

\begin{remark}
 Since $S\lsz n, d\rsz$ is clearly nonnegative, Lemma \ref{l1} implies
 that $k_d\geq 0$ for all $d\geq m$. This means that 
\[\sum_{j=0}^d c_j\geq \sum_{j=0}^d p_j \quad\lef d\geq m\ri.\]

Loosely speaking, the degree of a typical vertex is asymptotically
larger than or equal to the degree of the new vertex. This is in
accordance with the fact that the degree of a fixed vertex cannot
decrease. 

Similarly to Lemma 1, one can prove that 
\[\sum_{j=0}^d x_j\geq \sum_{j=0}^d q_j \ \ \ \lef d\geq m\ri,\]
which means the same for the selected vertices.\end{remark}


\section{Graph models}

In this section we briefly review some scale free random graph
models and sets of selected vertices to which the results of the
previous section can be applied.  

\subsection{Generalized plane oriented recursive tree}

We start from one edge, and at each step one new vertex and one new
edge are added to the graph. At the $n$th step the probability that a
given vertex of degree $d$ is connected to $v_n$ is
$\left(d+\beta\right)/T_{n-1}$, where $\beta>-1$ is the parameter of
the model, and $T_{n-1}=\left(2+\beta\right)(n+1)+\beta$. These kind of
random trees are widely examined, see for example
\cite{Drmota,Pittel}. $\beta=0$ gives the Albert--Barab\'asi tree
\cite{[Bar99]}.  

We fix an integer $j\geq 1$. At the $n$th step $v_n$ is added to the
set of selected vertices if it is at distance $j$ from $u_1$ in
$G_n$. Thus $S_n$ is the $j$th level of the tree $G_n$. 

It is well known \cite{Mo02} that Condition \ref{c1} is
satisfied with 
\[c_d=\frac{\lef 2+\beta\ri\Gamma\lef d+\beta\ri\Gamma\lef 3+2\beta\ri}
{\Gamma\lef 1+\beta\ri\Gamma\lef d+3+2\beta\ri }\quad\lef d\geq 1\ri.\]  
Consequently, 
\[c_d\sim \frac{\lef 2+\beta\ri\Gamma\lef 3+2\beta\ri}{\Gamma\lef
  1+\beta\ri}\cdot d^{-\lef 3+\beta\ri}\quad\lef d\tv\ri\] 
and $\gamma=3+\beta$ satisfies Condition \ref{c2}. It is clear
from the definition that Condition \ref{c3} holds, and since the
degree of the new vertex is always $1$, we have $m=1$, and conditions
\ref{c4}, \ref{c5}, and \ref{d4} are trivially
satisfied. Using that $p_d=0$ for $d\neq 1$ and $p_1=1$, Condition 
\ref{c6} is also easy to check. 

The distance of $v_n$ and $u_1$ does not change after generating the
edges from $v_n$ at the $n$th step. This guarantees Conditions
\ref{d1} and \ref{d2}. The results of \cite{[Mo06]} show
that Condition \ref{d3} is satisfied with
$\alpha=1/\left(2+\beta\right)$. It is proven that 
\[\left\vert S_n\right\vert \sim n\zeta\,\frac{\mu\lef n\ri^{j-1}}{\lef
j-1\ri!}\,e^{-\mu\lef n\ri}\asymp n^{\tfrac{1}{2+\beta}}\lef \log
n\ri^{j-1}\quad\lef n\tv\ri,\] 
where $\zeta$ is a positive random variable, and $\mu\lef
n\ri=\frac{1+\beta}{2+\beta}\log n$ \cite[Theorem 2.1]{[Mo06]}. 

Thus Theorem \ref{t1} applies: the asymptotic degree distribution
constrained on a fixed level of the tree does exist. The new
characteristic exponent is the following (cf. \cite[Theorem
3.1]{[Mo06]}).  
\[\gamma\cs=\al\lef \gamma-1\ri+1=\frac{1}{2+\beta}\lef 3+\beta-1\ri+1=2.\] 

\subsection{Independent edges}

We start from one edge. At the $n$th step, independently of each other,
every old vertex  is connected to the new one with probability
$\lambda d/T_{n-1}$, where $d$ is the degree of the old vertex in
$G_{n-1}$,  $0<\lambda<2$ is a fixed parameter, and $T_{n-1}$ denotes
the sum of degrees  in $G_{n-1}$. The restriction on $\lambda$
guarantees that the probability given above belongs to $\lsz
0,1\rsz$. It is clear that $m=0$. 

We fix one vertex, $v$, and $S_n$ consists of its neighbours in $G_n$. 

 In \cite[Theorem 3.1.]{[Ka06]} it is proven that 
the asymptotic degree distribution is given by
\[c_0=p_0,\quad c_d=\frac{2}{d\lef d+1\ri\lef d+2\ri}\sum_{k=1}^d k\lef k+1
\ri p_k,\]
where 
\[p_k=\frac{\lambda^k}{k!}\,e^{-\lambda}.\] 
Clearly, $c_d\sim 2\lambda\lef 2+\lambda\ri d^{-3} \quad\lef d\tv\ri$.
Thus the first two conditions are satisfied, and $\gamma=3$. Condition
\ref{c3} holds, because the probability that a given vertex gets
a new edge depends only on its actual degree. It is also clear that
Conditions \ref{d1} and \ref{d2} hold. Condition \ref{d3} is a
corollary of \cite[Theorem 2.1]{[Mo07]}, and we have $\alpha=1/2$. 

In this case the initial degree of the new vertex is not fixed. It is
proven in \cite{[Ka06]} that 
\[\sum_{k=0}^{\infty} \left\vert E\lef\left. I\lsz n+1, d\rsz\right\vert
\F_n\ri-p_d\right\vert\rightarrow 0\] 
almost surely as $n\tv$. 
This, and the fact that $\lef p_d\ri$ is a Poisson distribution with
parameter $\lambda$ imply Condition \ref{c4}.  

Note that the conditional distribution of $Y\lsz n, d\rsz$ is binomial
of order $X\lsz n, d\rsz$ and parameter $\lambda d/T_n<1$. One can
check Condition \ref{c5} with $Z_d$ having a suitable Poisson
distribution.

Condition \ref{d4} can be verified basing on the fact that the
degree distribution of a new selected vertex is similar to the
distribution of a new vertex because of the independent random
choices, and the following results. 
Theorem 2.1 in \cite{[Mo07a]} states that $T_n=2\lambda n+o\lef
n^{1-\e}\ri$ almost surely if $\e>0$ is sufficiently small. Moreover,
Theorem 2.2 there implies that the maximum degree after $n$ steps is
$O\lef \sqrt{n}\ri$ almost surely.

Our Theorem \ref{t1} can be applied, so the almost sure asymptotic
degree distribution constrained on the neighbours of a fixed vertex
exists. The new characteristic exponent is given by 
\[\gamma\cs=\alpha\lef \gamma-1\ri+1=\frac12\lef 3-1\ri+1=2\]
(cf. \cite[Theorem 3.1]{[Mo07]}).


Let us modify this example in such a way that vertices of degree $1$
never get new edges. Let
\[T_{n-1}=\sum_{d=2}^n X[n-1,d]\,d,\]
and choose $S_n$ to contain all vertices of degree $1$. Then we can see
that the all conditions hold except Condition \ref{c6}, but 
$x_d=0$ for $d>1$. This shows that positivity of $k_d$ cannot be
relaxed in order to obtain a polynomially decreasing degree distribution. 


\subsection{Random multitrees}

For $M\geq 2$ an $M$-multicherry is a hypergraph on $M+1$
vertices. One of them, called center, is distinguished, it is
connected to all other vertices with ordinary edges ($2$-hyperedges),
and the remaining $M$ vertices form an $M$-hyperedge, called the base. 

We start from the complete graph of $M$ vertices; the vertices form a
base. Then at each step we add a new vertex and an $M$-multicherry
with the new vertex in its center. We select the base of the new
multicherry from the existing bases uniformly. Finally, we add $M$ new
bases by replacing a vertex in the selected base with the new center in
all possible ways.  

The degree of the new vertex is always $M$, thus $m=M$.

Let $S_n$ be the set of vertices that are at distance $j$ from the
initial configuration. 

It is shown in \cite{[Mo09]} that Conditions \ref{c1},
\ref{c2}, and \ref{d3} are satisfied with
$\gamma=2+\frac{1}{M-1}$ and $\alpha=\frac{M-1}{M}$. The
other conditions are easy to check, using that distances in the
multitree do not change.  

Therefore Theorem \ref{t1} applies, and
\[\gamma\cs=\alpha\lef \gamma-1\ri+1=\frac{M-1}{M}\lef
2+\frac{1}{M-1}-1\ri+1=2.\] 

Another option for the set of selected vertices is the
following. Fix an integer $1\leq k<M$ and $k$ different
vertices. Let $S_n$ be the set of vertices that are connected to all
of them. Since the model is the same, we only have to check the
conditions on the set of selected vertices. Now Conditions \ref{d1},
\ref{d2}, and \ref{d4} clearly hold. Condition \ref{d3} can be proven
by slight modifications of the proofs of \cite{[Mo09]}. In this case
$\gamma=2+\frac1{M-1}$, $\alpha=1-\frac kM$, and 
\[\gamma\cs=2-\frac{k-1}{M-1}>1.\]

\section{Conclusions}
We presented sufficient conditions for the existence of the asymptotic
degree distribution constrained on the set of selected vertices. Scale
free property and regular variation of the size of the set of selected
vertices were essential. The new characteristic
exponent depended only on $\gamma$ and $\alpha$.  
 
We reviewed several models satisfying these conditions and identified
their characteristic exponents 
applying our main result. In these models $\gamma\cs\leq 2$ and
$\gamma\cs\leq \gamma$, thus the characteristic exponent decreased. One
reason for that is the following. The selected vertices are closer to
the initial configuration in some sense.  There  are more ``old''
vertices among them and their degree is larger than that of the
``typical'' ones.

\end{document}